\documentclass[12pt]{article}
\usepackage{amsfonts,amsmath,amssymb,amsthm}
\usepackage[francais,english]{babel}
\usepackage[T1]{fontenc}
\usepackage[latin1]{inputenc}  
\usepackage{color}

\usepackage[linesnumbered,vlined,ruled,english,onelanguage]{algorithm2e}
\usepackage{enumitem}
\usepackage{titlesec}
\usepackage{upgreek}

\usepackage{hyperref}
\usepackage{minitoc}
\usepackage{appendix}

\hypersetup{colorlinks, linkcolor=magenta}

\usepackage{fancyhdr}
\usepackage{float}
\usepackage{makeidx}
\usepackage[dvipsnames]{xcolor}

\usepackage{dsfont}

\usepackage{cite}     
\usepackage{xcolor}
\usepackage{graphics}
\usepackage{graphicx}

\usepackage{bigints}

\usepackage{caption}
\usepackage[active]{srcltx} 

\newtheorem{theorem}{Theorem}

\newtheorem{proposition}{Proposition}
\newtheorem{remark}{Remark}
\newtheorem{definition}{Definition}

\newcommand{\gras}{G_m(\mathbb{R}^n)} 
\newcommand{\egras}{\mathbb{R}^{n\times m}_*}
 
\DeclareMathOperator*{\argmin}{argmin}

\usepackage{array}
\newcolumntype{F}{ >{\centering\arraybackslash} m{3.5cm} } 
\newcolumntype{T}{ >{\centering\arraybackslash} m{1.5cm} } 

\usepackage{subcaption}

\newcommand{\vertiii}[1]{{\left\vert\kern-0.25ex\left\vert\kern-0.25ex\left\vert #1 
    \right\vert\kern-0.25ex\right\vert\kern-0.25ex\right\vert}}

\begin{document}

\title{\Large\sf Generalization of the Neville-Aitken Interpolation Algorithm on Grassmann Manifolds : Applications to Reduced Order Model}
\author{\sf 
R. Mosquera \footnote{Laboratoire LaSIE, University of La Rochelle, France. Email: rolando.mosquera$\_$meza@univ-lr.fr} , 
A. El Hamidi \footnote{Email: aelhamid@univ-lr.fr} , 
A. Hamdouni \footnote{Email: ahamdoun@univ-lr.fr} , 
A. Falaize \footnote{Email: antoine.falaize@univ-lr.fr}
}
\date{}
\maketitle

\renewcommand{\thefootnote}{\arabic{footnote}}

\begin{abstract}
The interpolation on Grassmann manifolds in the framework of 
parametric evolution partial differential equations is presented. Interpolation points on the Grassmann manifold 
are the subspaces spanned by the POD bases of the available solutions corresponding to the chosen parameter values. The well-known Neville-Aitken's algorithm is extended to Grassmann manifold, where interpolation is performed in a recursive way via the geodesic barycenter of two points. The performances of the proposed method are illustrated through three independent CFD applications, namely: the Von Karman vortex shedding street, the lid-driven cavity with inflow and the flow induced by a rotating solid. The obtained numerical simulations are pertinent both in terms of the accuracy of results and the time computation.
\end{abstract}




\pagestyle{myheadings}
\thispagestyle{plain}
\markboth{}
{On the interpolation on Grassmann Manifolds : Neville-Aitken's Algorithm}

\section{Introduction} \label{intro}
In the present manuscript, we are interested with interpolation on Grassmann manifolds in the framework of 
parametric evolution partial differential equations. A relevant example in this topic is the Navier-Stokes equations where the varying parameter $\lambda$ is the Reynolds number. 
%
%
We confine ourselves in the situation where the solution of the parametric partial differential equation is computed for a set $\{\lambda_0, \, \lambda_1, \cdots, \lambda_N\}$ of parameter values and the Proper Orthogonal Decomposition (POD) is used to construct the associated bases $\Phi(\lambda_0), \, \Phi(\lambda_1), \cdots \Phi(\lambda_N)$. 
These solutions are in general computed with high computational costs discretization type methods as Finite Elements, Finite Volumes, Discontinuous Galerkin, etc.

For a new value of the parameter $\lambda$, instead of recomputing the solution, we can use the already 
computed POD bases to carry out an interpolation. A natural question is how to interpolate a set of given bases. The purpose of the present work is to develop low-dimensional parametric methods to compute an approximation of the solution, by interpolation of POD bases, for a new parameter value.

In various problems in fluid mechanics or fluid structure interaction, the POD bases are
computed for a given Reynolds number or for other flow parameters. Theses bases are then
used to build a reduced model to predict the flow for other values of these parameters.
The parametric domain of validity of such reduced models has to be described.
In the framework of quasilinear parabolic problems, this parametric domain of validity was studied according to the solution's regularity and the number of modes retained to build the reduced model by Akkari {\em et al.} \cite{Akkari_AMSES}. 
These results were then extended to Navier-Stokes equations \cite{Akkari_al1} and Burgers equation \cite{Akkari_al2}.
Numerical results on parametric sensitivity for fluid structure interaction problems were carried out in \cite{Pomarede_thesis}. In \cite{Gravouil}, the authors introduced non-intrusive strategies to provide vademecums intended to real-time computations for nonlinear thermo-mechanical problems.

In control theory, the numerical costs (CPU
and memory) associated with the adjoint equation-based methods
used to solve the underlying optimization problems
are so important that the three-dimensional Navier-Stokes equations are rarely studied.

In \cite{Kunisch2}, the POD is applied to solve open-loop and closed-loop optimal control problems for the Burgers equation. Similar methods are used to control laser surface hardening, where the state equations are a semilinear heat equation coupled with an ordinary differential equation \cite{Homberg}. Phase-field control problems 
are also investigated by \cite{Volkwein}. In medical imaging \cite{BAC,Roujol}, parametric POD bases are studied 
when time-interval variation, spatial domain variation or parameter variation may occur.


One classical way to perform such approximations is the use of {\it Reduced Basis} (RB) methods. RB-methods consist in the approximation of the manifold solution $\mathcal{M} \subset H$, {\it i.e.,} the set of parametric solutions, by a low-dimensional subspace $H_{m} \subset H$, where $H$ is the Hilbert space in which the functional framework 
of the problem is well-posed. In this topic, one popular way is the construction of such a subspace via the {\it snapshots}
$\left( u(\lambda_k) \right)_{0 \leq k \leq {N}}$ as, for example, $H_{m} = {\rm span} \{u(\lambda_k) \; : \; 0 \leq k \leq {N} \}$ and
use Galerkin type techniques to compute an approximation $u_{_m} (\lambda) \in  H_{m}$ of
$u(\lambda) \in \mathcal{M}$. We refer the interested reader to \cite{Patera_Rozza,Haasdonk} and the references therein.

Another way to approximate the parametric solution $u$ for a new parameter value $\lambda \in \mathbb{R}^{P}$
from given data $\left( u(\lambda_k) \right)_{0 \leq k \leq {N}}$ via POD techniques is the interpolation
on {\it Grassmann manifolds}. We describe briefly this method and refer the reader to the work of Amsallem and Farhat~\cite{Amsallem2008} for more details. This procedure can be summarized as follows: for each $k \in \{0,1, \cdots, {N}\}$, consider the rank-$m$ POD basis $\left(\Phi_{_i} (\lambda_k) \right)_{1 \leq i \leq m}$ and the corresponding spanned subspace
$H_{m}^{(k)} = {\rm Span} \left \{\Phi_{_i} (\lambda_k) \; : \; 1 \leq i \leq m \right\} \subset H$,
(in \cite{Amsallem2008}, the space $H$ has a large but finite dimension). Now, denote by $G_m(H)$ the Grassmann manifold consisting of all $m$-dimensional subspaces of $H$. Then the $m$ subspaces
$H_{m}^{(k)}$ can be seen as points on this manifold $G_m(H)$. For a new parameter value
$\lambda$, the question is how to compute a good approximation of the rank-${m}$ POD basis
associated with the solution $u(\lambda)$, or in a more relevant way, how to find its corresponding point in the Grassmann manifold $G_m(H)$? Indeed, it is known from the work of Amsallem and Farhat~\cite{Amsallem2008} that the appropriate objects to interpolate are not the POD bases but the underlying spanned subspaces.
To answer this question, the idea developed in \cite{Amsallem2008} is to choose a reference parameter value
$\lambda_{r} \in \{\lambda_k \; : \; 0 \leq k \leq {N} \}$, to consider the POD subspace $R := H_{m}^{(r)}$ as a reference point on the Grassmann manifold $G_m(H)$. Let us denote by $T_{_R} G_m(H)$ the tangent space to the manifold $G_m(H)$ at the point $H_{m}^{(r)}$. Then, the logarithm mapping at the reference point $H_{m}^{(r)}$ denoted by
$\log^{(r)} \; : \; G_m(H) \longrightarrow T_{_R} G_m(H)$ allows to map (in a local way) all points
$\left( H_{m}^{(k)} \right)_{0 \leq k \leq N}$ on the tangent space $T_{_R} G_m(H)$ where one can perform an interpolation (see Section 2 for more details). Finally, one returns to the manifold $G_m(H)$ via the (reciprocal) exponential mapping from the tangent space $T_{_R} G_m(H)$ to get the desired approximation of $H_m^{(\lambda)}$ on $G_m(H)$.

The disadvantage of the method developed in \cite{Amsallem2008} is its dependence on the reference point and the lack of a strategy for choosing this reference point. \\

In a recent work \cite{Mosquera_DCDS18}, the authors developed the Grassmann Inverse
Distance Weighting (G-IDW) as an extension of the well-known Inverse Distance Weighting (IDW) method to Grassmann manifolds. 
It consists in the minimization of a quadratic function in the geodesic distance to interpolation points, with appropriate weights. This G-IDW method does not require a reference point and yields relevant results. As a counterpart, it is iterative. \\

In the present work, we propose the extension of the Neville-Aitken's algorithm which 
computes the Lagrange interpolation polynomial in a recursive way from the interpolation of 
two points. Replacing straight lines by geodesics, we extend this algorithm to Grassmann manifolds. 
The obtained numerical simulations are excellent both in terms of accuracy of the results and the computation time. \\

The present paper is organized as follows: In Section 2, the Grassmann manifold and the underlying tools used in the paper are briefly presented. Section 3 is dedicated to the introduction of the Neville-Aitken's algorithm in the framework of vector spaces and its extension to Grassmann manifolds. Finally, 
the performances of the proposed method are illustrated through three independent CFD applications, namely: the Von Karman vortex shedding street, the lid-driven cavity with inflow and the flow induced by a rotating solid. 
Our results are compared to those obtained by the exact POD method and to the interpolation approach 
developed by Amsallem et al. \cite{Amsallem2008}. 


\section{The Grassmann manifold}
%
Let $m$, $n$ be two positive integers with $m\leq n$. The set of all $\mathbb R$-subspaces of $\mathbb R^n$ of dimension $m$ is a differentiable manifold of dimension $m\times(n-m)$, called the Grassmann
manifold and usually denoted by $\gras$. Let $\mathbb{R}^{n\times m}_*$ be the set of all matrices of size
$n \times m $, whose column vectors are linearly independent; one defines on this set the following equivalence relation
$$
X \sim Y, \hspace{0.3cm}\text{if and only if}\hspace{0.3cm}\exists\: A \in GL_m(\mathbb{R})\:\text{such that :}\: X=YA,
$$
where $GL_m(\mathbb{R})$ denotes the linear group of degree $m$, {\it i.e.,} the set of $m \times m$ invertible matrices.
The Grassmann manifold can be realized as the set of such equivalence classes. More precisely,
$$
\gras= \{\:\overline{X}: \:\:X \in \egras\: \},
$$
where
$$
\overline{X} =\{X\:A: \:\: A\in GL_m(\mathbb{R})\}.
$$
The Grassmann manifold is a
compact topological space and, thanks to the canonical projection denoted by
$$
\pi:\egras\to \gras,
$$
it can be endowed by a differentiable manifold structure so that $\pi$ is a submersion (\cite{Milnor1974}).
%
\subsection{Exponential and logarithm mappings, injectivity radius}
%
On each point $X\in\egras$, one constructs a metric $g_{_X}$ which is right-invariant
by $GL_m(\mathbb{R})$, {\em i.e.,}
$$
g_{_{XA}}(Y_1A,Y_2A)=g_{_X}(Y_1,Y_2),
$$
defined by
\begin{equation} \label{metric}
g_{_X}(Y_1,Y_2)= tr [(^tXX)^{-1}\:\:^tY_1\:Y_2],
\end{equation}
for every $Y_1,Y_2\in \mathbb{R}^{n \times m}$ and $A\in GL_m(\mathbb{R})$. This metric induces a Riemannian metric on $\gras$ through the submersion $\pi$. Moreover, one constructs a Riemannian connection structure on $\gras$, which allows the definition of geodesics and the geodesic {\em exponential } mapping on $\gras$ \cite{Absil2004}.

In what follows, the distance on $\gras$, induced by the metric $g$, will be denoted by $d$, and 
the ball centered at $y$ of radius $r>0$ in $\gras$ will be denoted by $B_{_d} (y,r)$.

%

%
For every $X \in \egras$, we introduce $H_{_X} = \left[ T_{_X} \pi^{-1} \left(\overline{X}\right)\right]^{\perp}$, the orthogonal space 
to the tangent space at $X$ of the fiber $\pi^{-1} \left(\overline{X}\right)$, where the orthogonality is with respect to scalar product 
induced by the metric $g$. Then the mapping 
$$
d_X \pi \; : \; H_{_X} \longrightarrow T_{\overline{X}} \, \gras
$$
is an isomorphism. Now, given a vector $v \in T_{\overline{X}} \, \gras$, if we perform the singular value decomposition (SVD) of the matrix $(d_X \pi)^{-1}(v) \, \left(^tX X\right)^{1/2}$ as 
$$
(d_X \pi)^{-1}(v) \, \left(^tX X\right)^{1/2} = U \; \Sigma \; ^t V,
$$
the geodesic of initial point
$\gamma(0)=\overline{X}$ and initial velocity $\gamma^\prime(0)=v$ can be written in the form
$$
\gamma(s)=\overline{X(^tXX)^{-1/2}\:V\cos(\,s \Sigma)  + U\sin( \,s \Sigma)}.
$$

The geodesic exponential mapping is defined by $\exp_{\overline{X}}(v)=\gamma(1)$ and has
a matrix representation given by:
$$
X(^tXX)^{-1/2} V \cos (\Sigma) + U \sin(\Sigma).
$$
On the other hand, since $d_{0_{_{\overline{X}}}}\exp_{\overline{X}}=I$, then the inverse function theorem implies that $\exp_{\overline{X}}$ is a local diffeomorphism at $0_{_{\overline{X}}}$, the null vector of
$T_{\overline{X}} G_m \left(\mathbb{R} ^n \right)$.
In the same way, we can define the geodesic between two points
$\overline{X},\overline{Y}$ in $G_m(\mathbb{R}^n)$, with $^tX Y \in GL_m(\mathbb{R}) $, by
\begin{equation}\label{geo_points}
\gamma_{_{\bar{X},\bar{Y}}}(s)= \overline{X\:\widetilde{V} \cos\left(s\:\tan^{-1}\left(\widetilde{\Sigma}\right)\right) + \widetilde{U}\:\sin\left(s\:\tan^{-1}\left(\widetilde{\Sigma}\right)\right)}
,
\end{equation}

where $\widetilde{U}\:\widetilde{\Sigma}\:^t\widetilde{V}= Y(^tX\:Y)^{-1}-X$ is a SVD. \\
Notice that  
\begin{equation} \label{distance}
d(\overline{X},\overline{Y})=\|\tan^{-1}\left(\widetilde{\Sigma}\right)\|_{_F},
\end{equation}
where $\| \, \cdot \, \|_{_F}$ denotes the Frobenius norm on the space of matrices.
\begin{definition}[injectivity radius]
Let $\overline{X}\in \gras$. 
The injectivity radius of $(\gras,g)$ at $\overline{X}$ is defined by
$$
r_{_I} \left(\overline{X} \right)=\sup \left \{ \rho > 0 :\:\: \exp_{\overline{X}} :B(0_{\overline{X}},\rho)\to \exp_{\overline{X}}(B(0_{\overline{X}},\rho))\:\:\text{is a diffeomorphism} \right \},
$$
and the injectivity radius of $\gras$ is defined by
$$
r_{_I}=\inf\{r_{_I} \left(\overline{X} \right):\: \overline{X}\in \gras\}.
$$
\end{definition}
The value of the injectivity radius of $\gras$ is given by the following result.
\begin{theorem} \cite{Wong1967} \label{rayon}
Let $n, m$ be two positive integers such that $min\{m,n-m\} \geq 2$. Then
the injectivity radius of $\gras$ is given by $r_{_I} = \pi/2$.
\end{theorem}
The proof of Theorem \ref{rayon} can be found in \cite{Kozlov}. This Theorem gives the radius of the ball
where the reciprocal mapping of the geodesic exponential $\exp_{\overline{X}}^{-1}$ (or $\log_{\overline{X}}$) is defined. More precisely,  for $\overline{Y}\in\gras$, with $\displaystyle d(\overline{X},\overline{Y})<\frac{\pi}{2}$, then $\exp_{\overline{X}}^{-1}(\overline{Y})$ is represented by the matrix
$$
U_\alpha \tan^{-1}(\Sigma_\alpha)\: ^t V_\alpha,
$$
with the singular value decomposition
$$
U_\alpha \Sigma_\alpha ^t V_\alpha = [Y (^t X Y)^{-1}(^tXX)-X]\:(^tXX)^{-1/2}.
$$


In the next section, we will introduce an extension of the well-known Neville-Aitken algorithm \cite{Aitken,Neville} to the Grassmann manifold. 

%
%

\section{Neville-Aitken's algorithm for Lagrange interpolation}

In this section, we first recall the classical Neville-Aitken's algorithm for vector-valued functions, and second introduce our extension to the Grassmann manifolds.

\subsection{Recall of the classical Neville-Aitken's algorithm}
Consider a function $f \, : \, \Lambda \subset \mathbb{R} \longrightarrow V$,
where $V$ is a vector space, and a subset $\{\lambda_0, \lambda_1, \cdots\} \subset \Lambda$, with $f(\lambda_k) = y_k \in V$, for every $\lambda_k$.

Given an integer $N \geq 1$, the Lagrange 
interpolation polynomial of $f$ at $\{\lambda_0, \lambda_1, \cdots, \lambda_N\}$ will be denoted by 
$L[\lambda_0, \cdots, \lambda_N ; y_0, \cdots, y_N ]$. It is the unique polynomial of degree less or equal to $N$, satisfying:
\begin{equation} \label{LA}
L[\lambda_0, \cdots, \lambda_N ; y_0, \cdots, y_N ] (\lambda_k) = y_k, \;\; \text{for every} \;\; k \in \{0,1, \cdots, N\}.
\end{equation}
Let us consider the affine functions $\alpha_{_{[i,j]}}$ defined on $\Lambda$ by
\begin{equation} \label{alphaij}
\alpha_{_{[i,j]}} (\lambda) = \frac{\lambda - \lambda_i}{\lambda_j - \lambda_i}, \quad 
\forall \, i,  j \in \{0,1, \cdots\}.
\end{equation}
Then, the Neville-Aitken's algorithm can be summarized by the following induction formula:
\begin{eqnarray} \label{NA}
 L[\lambda_0, \cdots, \lambda_N ; y_0, \cdots, y_N] (\lambda) = 
&(1-\alpha_{_{[0,N]}}(\lambda))&  \, L[\lambda_0, \cdots, \lambda_{N-1} ; y_0, \cdots, y_{N-1}] (\lambda) + \\
&+ \; \alpha_{_{[0,N]}} (\lambda)& \, L[\lambda_1, \cdots, \lambda_N; y_1, \cdots, y_N] (\lambda),
\end{eqnarray}
with the initialization
$$
L[\lambda_k; y_k] (\lambda) = y_k \;\;  \text{for every} \; \; k \in \{0,1, \cdots, N\}.
$$
That is, for a fixed $\lambda \in \Lambda$, the quantity $L[\lambda_0, \cdots, \lambda_N ; y_0, \cdots, y_N] (\lambda)$ is the barycenter of \\
$L[\lambda_0, \cdots, \lambda_{N-1} ; y_0, \cdots, y_{N-1}] (\lambda)$ and $L[\lambda_1, \cdots, \lambda_N ; y_1, \cdots, y_N] (\lambda)$, weighted by
$(1-\alpha_{_{[0,N]}}(\lambda))$ and $\alpha_{_{[0,N]}} (\lambda)$,
respectively. 
Therefore, Neville-Aitken's algorithm reduces the Lagrange interpolation into recursive 
barycenter computations of two points at each step, with appropriate weights. Thus, for a fixed value of the parameter $\lambda$, we get
\begin{eqnarray*}
L[\lambda_0, \lambda_1 ; y_0, y_1] (\lambda) 
& = & 
(1-\alpha_{_{[0,1]}}(\lambda)) \, y_0 + 
\alpha_{_{[0,1]}} (\lambda) \, y_1, \\ 
L[\lambda_1, \lambda_2 ; y_1, y_2] (\lambda)
& = & 
(1-\alpha_{_{[1,2]}}(\lambda)) \, y_1 + 
\alpha_{_{[1,2]}} (\lambda) \, y_2,  \\
L[\lambda_0,\lambda_1, \lambda_2 ; y_0, y_1, y_2] (\lambda)
& = & 
(1-\alpha_{_{[0,2]}}(\lambda)) \, L[\lambda_0, \lambda_1 ; y_0, y_1] (\lambda) + 
\alpha_{_{[0,2]}} (\lambda) \, L[\lambda_1, \lambda_2 ; y_1, y_2] (\lambda),
\end{eqnarray*}
with ${\rm deg} (L[\lambda_0, \lambda_1 ; y_0, y_1])\leq 1,$ 
${\rm deg} (L[\lambda_1, \lambda_2 ; y_1, y_2])\leq 1,$
${\rm deg} (L[\lambda_0,\lambda_1, \lambda_2 ; y_0, y_1, y_2])\leq 2,$
and so on to obtain the Lagrange interpolation polynomial $L[\lambda_0, \cdots, \lambda_N ; y_0, \cdots, y_N]$, whose degree is less or equal to $N$.
\begin{remark}
\label{rem:neville_vector}
For vector-valued functions, the Lagrange interpolation polynomial is invariant by 
permutations of interpolation points. That is, if $\sigma$ is a permutation of the set $\{\lambda_0, \lambda_1, \cdots, \lambda_N\}$, then
$L[\lambda_{\sigma{(0)}}, \cdots, \lambda_{\sigma{(N)}} ; y_{\sigma{(0)}}, \cdots, y_{\sigma{(N)}}] = L[\lambda_0, \cdots, \lambda_N ; y_0, \cdots, y_N]$.
\end{remark}
\subsection{Extension of the Neville-Aitken's algorithm to Grassmann manifolds} \label{Neville}
As pointed before, if we are able to define the notion of barycenter of two points on the Grassmann manifold, 
then we will be able to carry out interpolation on it, via Neville-Aitken's algorithm. 

More precisely, consider a function $f \, : \, \Lambda \subset \mathbb{R} \longrightarrow \gras$ such that 
$y_i := f(\lambda_i)$ are given, for $i \in \{0,1, \cdots \}$. We will describe one way to handle with 
the interpolation of $f$ at $\{\lambda_0, \lambda_1, \cdots, \lambda_N\}$, for any $N \geq 1$.
\begin{remark}
Since our interpolation extension will be based on geodesics on the Grassmann manifold $\gras$, 
we will assume hereafter that all the points $\left( y_i\right)_{0 \leq i \leq N}$ belong to 
a certain ball of radius $\frac{\pi}{2}$ in $\gras$. That is 
\begin{equation} \label{hypo}
\text{({\bf H})} \quad \exists y \in \gras \; : \; \forall k \in \{0, 1, \cdots, N\}, \; 
y_k \in B_d(y,\pi/2).
\end{equation}
The assumption ({\bf H}) will ensure the uniqueness of 
geodesics connecting each pair of interpolation points \cite{Wong1967}.
\end{remark}
At first, consider two (distinct) values $\lambda_i$ and $\lambda_j$ whose images $y_i$ and $y_j$ belong to the Grassmann manifold 
$\gras$. Let $\gamma_{_{[y_i,y_j]}} \, : \, [0,1] \longrightarrow \gras$ 
be the unique geodesic on $\gras$ satisfying $\gamma_{_{[y_i,y_j]}}(0) = y_i$ and $\gamma_{_{[y_i,y_j]}}(1) = y_j$. Then, the map defined on $[\lambda_i,\lambda_j]$ by 
\begin{equation} \label{Gamma}
\Gamma [\lambda_i,\lambda_j \, ; \,  y_i, y_j] := \gamma_{_{[y_i,y_j]}} \circ \alpha_{_{[i,j]}}
\end{equation}
is the unique geodesic on 
$\gras$ satisfying $\Gamma [\lambda_i,\lambda_j \, ; \,  y_i, y_j](\lambda_i) = y_i$ and $\Gamma [\lambda_i,\lambda_j \, ; \,  y_i, y_j](\lambda_j) = y_j$, where $\alpha_{_{[i,j]}}$ is defined by (\ref{alphaij}). Indeed, geodesics 
are invariant by affine transformations
and thanks to equation (\ref{geo_points}) a matrix representation of the geodesic connecting two points on the manifold is known. This matrix representation will be used in the algorithm~\ref{alg_neville}.
\begin{definition} \label{gen_bar}
For a given $\lambda \in [\lambda_i, \lambda_j]$, the point 
$\Gamma [\lambda_i,\lambda_j \, ; \,  y_i, y_j] (\lambda) \in \gras$ is called the barycenter (in the sense of geodesics) on $\gras$ of \, $y_i$ and $y_j$, 
weighted by $\left(1 - \alpha_{_{[i,j]}} (\lambda)\right)$ and $\alpha_{_{[i,j]}} (\lambda)$ respectively. 
\end{definition}
The barycenter $\Gamma [\lambda_i,\lambda_j \, ; \,  y_i, y_j]$ provides an interpolation on $\gras$ of 
the function $f$ at $\lambda_i$ and $\lambda_j$. Moreover, it coincides with
 the barycenter in the sense of Karcher \cite{Karcher}. More precisely, we have the following result.
 %
 %
\begin{proposition}
Let $\lambda_i < \lambda_j$ be two parameter values in $\{\lambda_0, \lambda_1, \cdots  \}$ and 
$y_i$, $y_j$ be two points on the Grassmann manifold $\gras$.
Then
$$
\Gamma [\lambda_i, \lambda_j \, ; \,  y_i, y_j] (\lambda) = \underset{y \in \gras}{\arg\min} J(y),
$$
where 
\begin{equation} \label{J_karcher}
J(y) = \frac{1}{2} \left[ \left(1 - \alpha_{_{[i,j]}} (\lambda)\right) \, d^2 (y , y_i) + \alpha_{_{[i,j]}} (\lambda) \, d^2 (y , y_j) \right],
\end{equation}
and $d$ denotes the geodesic distance on $\gras$.
\end{proposition}
{\bf Proof.} Applying  \cite[Theorem 4.1]{Mosquera_DCDS18}, it suffices to show that $\nabla J (\Gamma [\lambda_i, \lambda_j \, ; \,  y_i, y_j] (\lambda))=0$, or equivalently:
$$
\left(1 - \alpha_{_{[i,j]}} (\lambda) \right) \exp_{_{\Gamma [\lambda_i, \lambda_j \, ; \,  y_i, y_j] (\lambda)}}^{-1}(y_i) +  \alpha_{_{[i,j]}} (\lambda) \exp_{_{\Gamma [\lambda_i, \lambda_j \, ; \,  y_i, y_j] (\lambda)}}^{-1}(y_j)=0.
$$
To lighten the notations, we set $\Gamma [\lambda_i, \lambda_j \, ; \,  y_i, y_j] (\lambda)= \widehat{y}$. 
We introduce the two geodesic paths on $\gamma_1$, and $\gamma_2$ defined from the interval $[0,1]$ to 
$\gras$ by 
\begin{eqnarray*}
\gamma_1 (t) &=& \Gamma [\lambda_i, \lambda_j \, ; \,  y_i, y_j] (\lambda + t (\lambda_i - \lambda)) \\
\gamma_2 (t) &=& \Gamma [\lambda_i, \lambda_j \, ; \,  y_i, y_j] (\lambda + t (\lambda_j - \lambda))
\end{eqnarray*}
from $\widehat y$ to $y_i$ and from $\widehat y$ to $y_j$, respectively. Therefore, we get
\begin{eqnarray}
\gamma_1'(0) &=& \exp_{\widehat{y}}^{-1} (y_i) = 
(\lambda_i - \lambda) \Gamma' [\lambda_i, \lambda_j \, ; \,  y_i, y_j] (\lambda) \label{gamma1}\\
\gamma_2'(0) &=& \exp_{\widehat{y}}^{-1} (y_j) = 
(\lambda_j - \lambda) \Gamma' [\lambda_i, \lambda_j \, ; \,  y_i, y_j] (\lambda). \label{gamma2}
\end{eqnarray}
Whence, combining (\ref{gamma1}) and (\ref{gamma2}), we obtain
$$
\left(1 - \alpha_{_{[i,j]}} (\lambda) \right) \exp_{\widehat{y}}^{-1}(y_i) +  \alpha_{_{[i,j]}} (\lambda) \exp_{\widehat{y}}^{-1}(y_j) = 0,
$$
which ends the proof. \hfill $\Box$ \\

In the same way, we can extend the analogous of (\ref{NA}) to determine the interpolation $\Gamma [\lambda_i,\lambda_j,\lambda_k \, ; \, y_i, y_j, y_k]$ of $f$ at $\{\lambda_i,\lambda_j,\lambda_k\}$. 
Indeed, we compute at first $\Gamma [\lambda_i,\lambda_j \, ; \, y_i, y_j] (\lambda)$ and $\Gamma [\lambda_j,\lambda_k \, ; \, y_j, y_k] (\lambda)$ by (\ref{Gamma}) and then 
\begin{equation*}
\Gamma [\lambda_i,\lambda_j,\lambda_k \, ; \, y_i, y_j, y_k] (\lambda) = 
\Gamma \big[\lambda_i,\lambda_k \, ; \, \Gamma [\lambda_i,\lambda_j \, ; \, y_i, y_j] (\lambda) , 
\Gamma [\lambda_j,\lambda_k \, ; \, y_j, y_k] (\lambda) \big] (\lambda),
\end{equation*} 
and so on to obtain the interpolation $\Gamma [\lambda_0,\cdots,\lambda_N \, ; \, y_0, \cdots, y_N]$ of $f$ at $\{\lambda_0, \cdots, \lambda_N\}$, that is :
\begin{small}
\begin{equation*}
\Gamma [\lambda_0,\cdots,\lambda_N \, ; \, y_0, \cdots, y_N] (\lambda) = 
\Gamma \big[\lambda_0,\lambda_{N} \, ; \, 
\Gamma [\lambda_0, \cdots, \lambda_{N-1} \, ; \, y_0, \cdots y_{N-1}] (\lambda) , 
\Gamma [\lambda_1, \cdots, \lambda_{N} \, ; \, y_1, \cdots y_{N}] (\lambda) \big] (\lambda).
\end{equation*} 
\end{small}
Therefore, the resulting Neville-Aitken algorithm on Grassmann manifolds can be summarized as follows: \\

\begin{algorithm}[H]
\DontPrintSemicolon
\KwData{
\begin{enumerate}
\item[]{$i.$} List of $(N+1)$ parameter values :
$\Lambda=\{\lambda_0, \lambda_1 \, \cdots,\lambda_N\} \, \subset \mathbb{R}$.
\item[]{$ii.$} List of $(N+1)$ bases: 
$\{\Phi(\lambda_0), \Phi(\lambda_1) \, \cdots,\Phi(\lambda_N)\} \, \subset \mathbb{R}^{n\times m}_*$.
\item[]{$iii.$} The new parameter value:
$\lambda\in \displaystyle \left\{\min_{0 \leq k \leq N} \lambda_k \, , \, \max_{0 \leq k \leq N} \lambda_k \right\}$.

\end{enumerate}
}

\vspace*{0.25cm}

\KwResult{The interpolated base $\Phi^I(\lambda)$}

\vspace*{0.25cm}

\DontPrintSemicolon
\For{$ i = 0$ to $N$}{
$Y (i,0) \longleftarrow \Phi(\lambda_i)$;
}
\For{$ j = 0$ to $N - 1$}{
\For{$ i = 0$ to $N - j - 1$}{
$Y(i,j+1) \longleftarrow \Gamma [\lambda_i,\cdots,\lambda_{i+j+1} \, ; \, y_i, \cdots, y_{i+j+1}] (\lambda)$
}}
$\Phi^I(\lambda) \longleftarrow  Y(0,N) $
\vspace*{0.25cm}
\caption{Neville-Aitken's algorithm on Grassmann manifolds}
\label{alg_neville} 
\end{algorithm}

\vspace*{0.25cm}

%
%
\section{Applications}
The performances of the proposed method are illustrated in this section through three CFD applications, namely (i) the Von Karman vortex shedding street, (ii) the lid-driven cavity with inflow and (iii) the flow induced by a rotating solid.
For each application, a relevant dimensionless parameter $\lambda$ is selected and the test procedure is as follows.
\begin{enumerate}
    \item We derive a High-Dimensional Model (HDM) by a standard finite-element discretization of the governing equations.
    \item We build $N_\lambda$ POD bases by applying the \emph{method of snapshots} \cite{Sirovich1987} to the solution arising from the simulation of the HDM for a set of $N_\lambda$ sampling parameters $(\lambda_i)_{i=0}^{N_\lambda}$.
    \item We interpolate the sets of POD bases at a given new parameter $\lambda^\star$ by (i) the proposed method (labeled \textsc{Neville}) and (ii) that proposed in~\cite{Amsallem2008} (labeled \textsc{Amsallem}). In the latter case, the reference point of the method $\Phi(\lambda_r)$ is chosen as the point in the sample bases closest to the target parameter: $\lambda_r=\lambda_i$ with $i = \underset{i \in \{0, \cdots, N_\lambda\}}{\argmin} \vert \lambda^\star - \lambda_i\vert$.
    \item We build a Reduced Order Model (POD-ROM) by a Gelerkin projection of the governing equations on the POD bases obtained from both interpolation methods.
    \item The results are compared with the POD-ROM built from the \emph{reference} POD basis computed directly from the snapshots of the HDM for $\lambda ^\star$ (labeled \textsc{Reference}).
\end{enumerate}
Details on the POD can be found in the reference~\cite{Holmes1997} and in~\cite{Volkwein2013, Cordier2003}.
All the numerical tests have been performed using the Python/C++ finite element library \textsc{DOLFIN}~\cite{logg2012dolfin} on a computer equipped  with a processor\footnote{2 sockets, 8 cores for each socket, 2 threads for each core, cadenced at 2.10GHz with a cache of 20MB.} Intel Xeon E5-2620 v4 and 64Go of RAM.
\subsection{Von Karman vortex shedding street\label{sec:canal}}
\begin{figure}
\includegraphics[width=\linewidth]{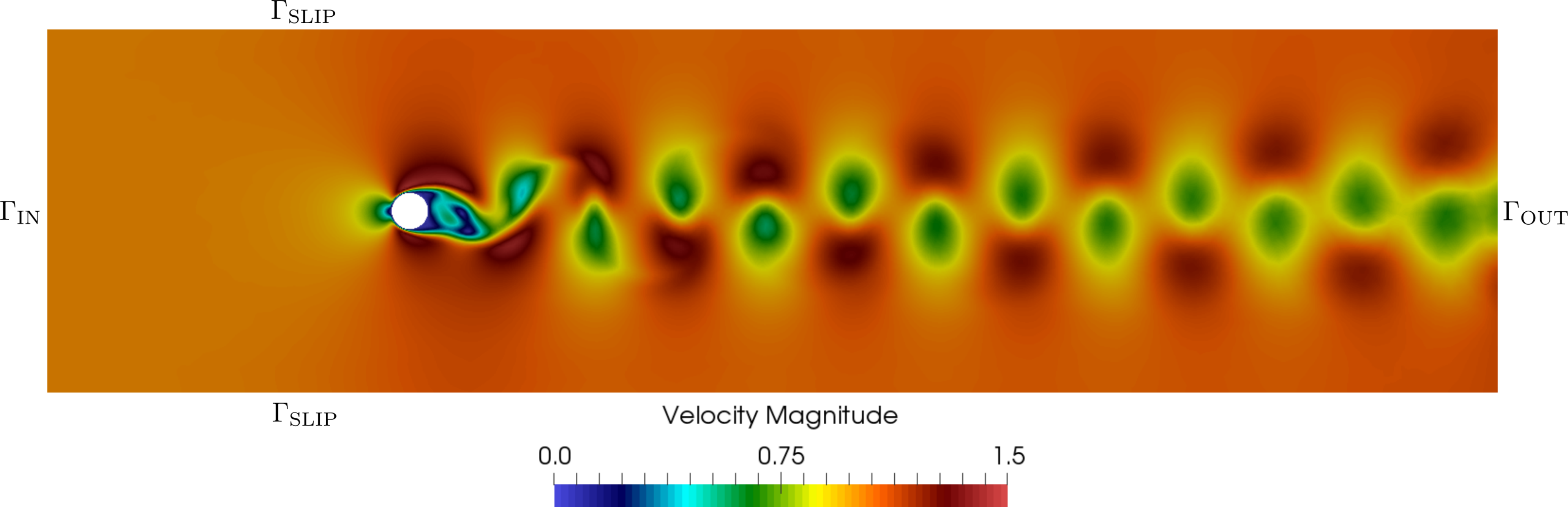}
\caption{\label{fig:canal} Description of the configuration for the application of section~\ref{sec:canal} (Von Karman vortex shedding street). The horizontal (respectively vertical) dimension of the channel is 40m (respectively 10m).
The solid disk is centered at $\mathbf{x} = (10, 5)$ and has a diameter $D=1$m.
The boundary conditions are inflow $\mathbf{v}=(1, 0)$ on $\Gamma_{\mathrm{IN}}$, outflow $\boldsymbol{\sigma} \cdot \mathbf{n}=0$ on $\Gamma_{\mathrm{OUT}}$ (with $\boldsymbol{\sigma}$ the Cauchy stress tensor and $\mathbf{n}$ the outward normal unit vector on the boundary), slip condition $\mathbf{v}\cdot \mathbf{n}=0$ on $\Gamma_{\mathrm{SLIP}}$ and no-slip condition $\mathbf{v} = (0, 0)$ on the disk boundary.
}
\end{figure}

Here, we consider the planar flow of an incompressible newtonian fluid around a solid disk in a channel.
The configuration is shown in figure~\ref{fig:canal}.
The HDM is derived by a standard variational formulation of the adimensional Navier-Stokes equations over the Taylor-Hood finite-elements space ($\mathrm{P2}$ for velocity and $\mathrm{P1}$ for the pressure).
Here, the dimensionless parameter is the Reynolds number $\lambda = {\frac{\rho \, v_{\infty} \, \ell}{\eta}}$ where $\rho=1$kg.m$^{-3}$ is the fluid density, $v_{\infty}=1$m.s$^{-1}$ is the magnitude of the velocity at inflow, $\ell=1$m is the disk diameter and $\eta$ is the dynamic viscosity.
The test cases (\textit{i.e.} the set of sampling parameters and the target parameter) are given in table~\ref{tab:cases_vonkarman}.
\begin{table}[h!]
	\centering
	\begin{tabular}{ccc}
		\hline
		Case & Sampling $\Lambda$ & Target $\lambda^\star$ \\
		\hline
		$1$ & $(100, 120, 130, 160, 170, 200)$ & $110$ \\ 
		$2$ & $(100, 160, 170, 180, 200)$ & $110$ \\ 
		$3$ & $(100, 120, 130, 140, 200)$ & $190$ \\
		\hline	
	\end{tabular}
	\caption{\label{tab:cases_vonkarman} 
	Description of the cases for the application of section~\ref{sec:canal} (Von Karman vortex shedding street). The set of sampling parameters for which the simulation of the full order model is performed is $\Lambda$. The interpolation is performed at the new parameter $\lambda^\star$.}
\end{table}

For each sampling parameter $\lambda \in \Lambda$, the simulation is performed by a standard predictor/corrector scheme (see \textit{e.g.} \cite[\textsection{7.3}]{ferziger2012computational}) and a set of ${N_T=200}$ regularly spaced snapshots of both the velocity $\big(\mathbf{v}(\mathbf{x}, t_i, \lambda)\big)_{i=1}^{N_T}$ and the pressure $\big(p(\mathbf{x}, t_i, \lambda)\big)_{i=1}^{N_T}$ are selected over a period $T={20}$s which {does not} include the transient period.
To ensure homogeneous boundary conditions for the POD modes, both the velocity and the pressure are split into a mean field $\overline{\mathbf{v}}(\mathbf x, \lambda) = \frac{1}{N_T}\sum_{i=1}^{N_T}\mathbf{v}(\mathbf{x}, t_i, \lambda)$ and a fluctuating field:
\begin{eqnarray}
\widetilde{\mathbf{v}}(\mathbf{x}, t_i, \lambda) &=&  {\mathbf{v}}(\mathbf x, t_i, \lambda) - \overline{\mathbf{v}}(\mathbf x, \lambda), \\
\widetilde{{p}}(\mathbf{x}, t_i, \lambda) &=& {{p}}(\mathbf x, t_i, \lambda) - \overline{{p}}(\mathbf x, \lambda).
\end{eqnarray}
Then, we compute the POD basis $\Phi(\mathbf x, \lambda)=\big(\phi_m(\mathbf x, \lambda)\big)_{m=1}^{N_T}$  (respectively $\Psi(\mathbf x, \lambda)=\big(\psi_m(\mathbf x, \lambda)\big)_{m=1}^{N_T}$) for the fluctuating velocity (respectively pressure) associated with each parameter by the method known as snapshots-POD \cite{Sirovich1987}.
The POD bases are truncated to $M$ modes for every parameters $\lambda \in \Lambda$ so that the instantaneous fields are approximated as 
\begin{eqnarray}
\widehat{\mathbf{v}}(\mathbf x, t, \lambda) = \overline{\mathbf{v}}(\mathbf x, \lambda) + \sum_{m=1}^{M}a_m(t, \lambda)\,\phi_m(\mathbf{x}, \lambda), \label{eq:ersatz_velocity}\\
\widehat{{p}}(\mathbf x, t, \lambda) = \overline{{p}}(\mathbf x, \lambda) + \sum_{m=1}^{M}b_m(t, \lambda)\,\psi_m(\mathbf{x}, \lambda). \label{eq:ersatz_pressure}
\end{eqnarray}

In this applcation, $M=10$.
We use three methods to interpolate the sets $\big(\Phi(\mathbf x, \lambda )\big)_{\lambda \in \Lambda}$ and $\big(\Psi(\mathbf x, \lambda )\big)_{\lambda \in \Lambda}$ at the parameter $\lambda^\star$: (i)~the grassmannian Neville-Aitken interpolation method proposed in subsection~\ref{Neville} (labeled \textsc{Neville}), (ii)~the method proposed by Amsallem and Farhat in~\cite{Amsallem2008} (labeled \textsc{Amsallem}) and (iii)~ a naive piecewise affine interpolation of the components associated with the sampled POD bases seen as regular matrices (labeled \textsc{Standard}).
We firstly examine the projection error for the fluctuating velocity defined for a given POD basis $\Phi$ as follows:
\begin{equation}
\label{eq:projection_error_definition}
\varepsilon_{\Phi}^{\mathrm{proj}} \triangleq 
\frac{
	\left\Vert 
	\widetilde{\mathbf V}_{\lambda ^\star} - P_\Phi \widetilde{\mathbf V}_{\lambda ^\star}
	\right\Vert^2_F
}{
\left\Vert
\widetilde{\mathbf V}_{\lambda ^\star}
	\right\Vert^2_F
},
\end{equation}
where the column of the matrix $\widetilde{\mathbf V}_{\lambda}$ are the snapshots of the discrete fluctuating velocity field obtained by the simulation of the HDM for the parameter $\lambda$, $P_\Phi \widetilde{\mathbf V} = \sum_{m=1}^M (\phi_m\vert\widetilde{\mathbf V})_{L^2(\Omega)}\,\phi_m$ denotes the matrix whose columns are the columns of $\widetilde{\mathbf V}$ projected onto the subspace engendered by the POD basis $\Phi$, and $\left\Vert \, \cdot \, \right\Vert_F$ denotes the Frobenius norm.
The results are shown in table~\ref{tab:reconstruction_error_vonkarman} where the benefit of both interpolations based on the Grassmann manifold is evident when compared with the standard interpolation of the bases components. 
Also, the proposed grassmannian Neville-Aitken interpolation yields globally a lower error than the method proposed in~\cite{Amsallem2008}.
\begin{table}[h!]
	\centering
	\begin{tabular}{llll}
		\hline
		Method & Case 1 & Case 2 & Case 3 \\
		\hline
		\textsc{Reference} &  6.89e-05  &  6.89e-05  &  3.10e-04 \\
		\textsc{Neville} & { 1.03e-03 } & { 7.95e-03 } & { 1.64e-03 } \\
		\textsc{Amsallem} & { 2.54e-03 } & { 9.12e-02 } & { 6.68e-03 } \\
		\textsc{Standard} & { 8.51e-01 } & { 2.13e-01 } & { 1.41e-01 } \\
		\hline
	\end{tabular}
	\caption{
	\label{tab:reconstruction_error_vonkarman}
	Results for the projection error defined in (\ref{eq:projection_error_definition}) for the application of section~\ref{sec:canal} (Von Karman vortex shedding street).
	\textsc{Reference} refers to the basis built directly from the snapshots obtained by the simulation of the HDM at the target parameter $\lambda ^\star$ (no interpolation). \textsc{Neville}, \textsc{Amsallem} and \textsc{Standard} refer to the bases obtained from the interpolation of the sampling sets defined in table~\ref{tab:cases_vonkarman} (no simulation of the HDM).}
\end{table}
Then we construct so-called POD-ROMs by a standard Galerkin projection of the momentum equation onto (i) the POD basis for the velocity and (ii) the gradient of the POD basis for the pressure (see \textit{e.g.} \cite{Tallet2015}), with the bases labeled \textsc{Reference} (direct POD of the HDM solution for the target parameter), \textsc{Neville}, \textsc{Amsallem} and \textsc{Standard}.
This yields dynamical systems for the temporal coefficients $(a_m)_{m=1}^{M}$ and $(b_m)_{m=1}^{M}$ in (\ref{eq:ersatz_velocity}--\ref{eq:ersatz_pressure}).
We examine the error associated with the simulation of these dynamical systems defined as follows:
\begin{equation}
\label{eq:dynamic_error_definition}
\varepsilon_{\Phi}^{\mathrm{dyn}} \triangleq 
\frac{
	\left\Vert 
	\widetilde{\mathbf V}_{\lambda ^\star} - \widehat{\mathbf V}_{\Phi}
	\right\Vert^2_F
}{
\left\Vert
\widetilde{\mathbf V}_{\lambda ^\star}
	\right\Vert^2_F
},
\end{equation}
where the matrix $\widehat{\mathbf V}_{\Phi}$ is the reconstruction of the fluctuating velocity $\widehat{\mathbf v} - \overline{\mathbf v}$ in~(\ref{eq:ersatz_velocity}) for the coefficients $(a_m)_{m=1}^{M}$ obtained by the simulation of the POD-ROM associated with the basis $\Phi$.
The results are shown in table \ref{tab:dynamic_error_vonkarman} where we see that the POD-ROMs built from bases interpolated by the proposed grassmannian Neville-Aitken interpolator yields better results compared with the method proposed in~\cite{Amsallem2008}.
\begin{table}[h!]
	\centering
	\begin{tabular}{llll}
		\hline
		Method & Case 1 & Case 2 & Case 3 \\
		\hline
		\textsc{Reference} &  4.41e-04  &  4.41e-04  &  1.37e-03 \\
		\textsc{Neville} & { 1.74e-03 } & { 1.35e-02 } & { 4.46e-03 } \\
		\textsc{Amsallem} & { 3.42e-03 } & { 1.05e-01 } & { 1.31e-02 } \\
		\textsc{Standard} & { 9.36e-01 } & { 3.06e-01 } & { 1.60e-01 }  \\
		\hline
	\end{tabular}
	\caption{
	\label{tab:dynamic_error_vonkarman}
	Results for the error associated with the dynamical systems as defined in (\ref{eq:dynamic_error_definition}) for the application of section~\ref{sec:canal} (Von Karman vortex shedding street).}
\end{table}

\begin{figure}[h!]
\includegraphics[width=\linewidth]{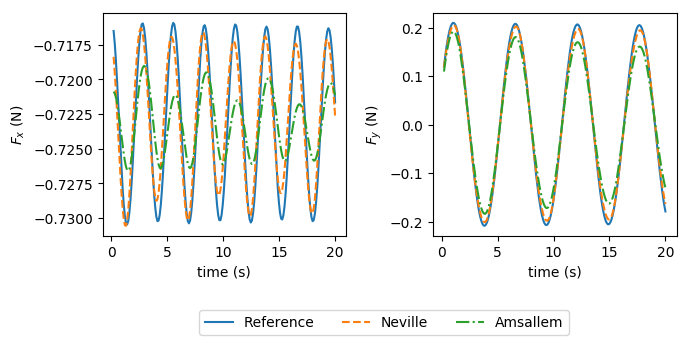}
\caption{
		\label{fig:efforts_aero_cylindre} 
		Drag effort $F_x$ and lift effort $F_y$ on the disk boundary for the Von Karman vortex shedding street of section~\ref{sec:canal}. 
}
\end{figure}

Finally, we compare the aerodynamic efforts imposed by the fluid on the body, defined as follows:
\begin{equation}
\label{eq:def_efforts_aero}
\left(\begin{array}{c}
F_x \\
F_y
\end{array}
\right) = \int_{\Gamma_{\mathrm{disk}}} \sigma(\mathbf{x}, t) \cdot \mathbf{n}(\mathbf{x}) \, \mathrm{d}\ell,
\end{equation}
where $\Gamma_{\mathrm{disk}}$ denotes the disk boundary, $\sigma$ is the fluid stress tensor, $\mathbf{n}$ is the outward unit vector normal to $\Gamma_{\mathrm{disk}}$, and $\mathrm{d}\ell$ is a line element.
Results are shown in figure~\ref{fig:efforts_aero_cylindre} for the case 2, where we see that the efforts obtained by the POD-ROM built from bases interpolated by the proposed grassmannian Neville-Aitken interpolator are closer to the reference.

\subsection{Lid-driven cavity with inflow\label{sec:liddriven}}
\begin{figure}[h!]
    \centering
    \begin{subfigure}{.35\linewidth}
        \includegraphics[width=\textwidth]{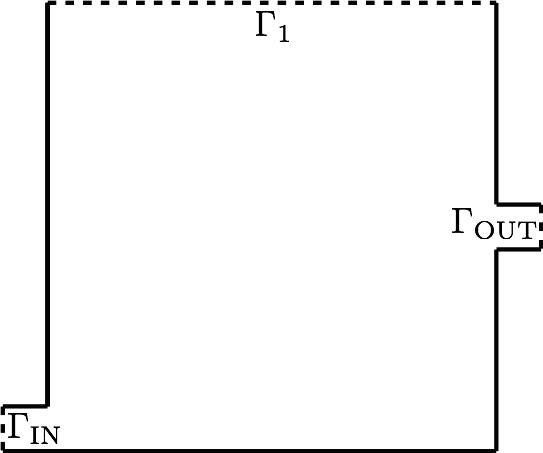}      
        \caption{\label{fig:lidriven_schematic} Schematic.}
    \end{subfigure}
    \begin{subfigure}{.35\linewidth}
        \includegraphics[width=\textwidth]{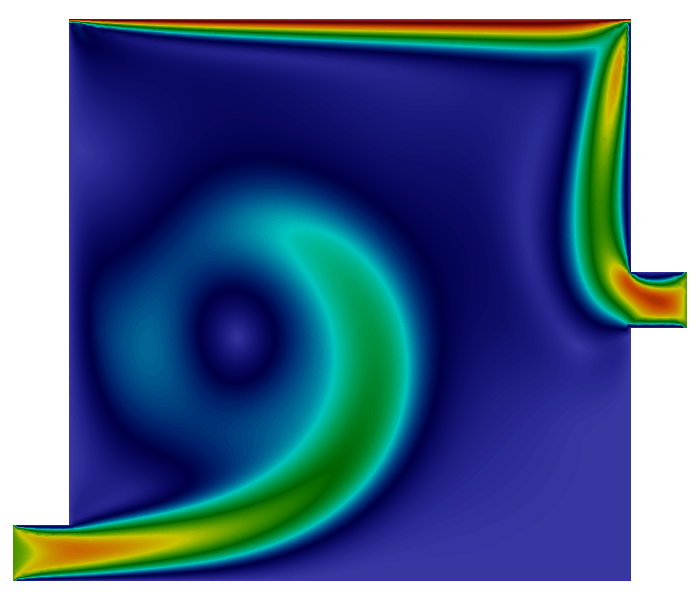}      
        \caption{\label{fig:lidriven_snapshot} Snapshot}
    \end{subfigure}
    \caption{\label{fig:liddriven_configuration} 
Description of the configuration for the application of section~\ref{sec:liddriven} (Lid-driven cavity with inflow). 
The dimensions of the cavity are 1m$\times$1m and the inflow/outflow extends by 0.1m.
The imposed velocity is $\mathbf{v}\vert_{\Gamma_1}=(1, 0)$ and $\mathbf{v}\vert_{{\mathrm{IN}}} = (v_x, 0)$ with $v_x \equiv \lambda$ and $\lambda \in (0, 1)$. 
The outflow boundary condition is $\boldsymbol{\sigma} \cdot \mathbf{n}=0$ on $\Gamma_{\mathrm{OUT}}$ (with $\boldsymbol{\sigma}$ the Cauchy stress tensor and $\mathbf{n}$ the outward normal unit vector on the boundary).
Homogeneous Dirichlet boundary conditions are imposed on the remaining boundary $ \Gamma_{\mathrm{0}} =  \partial \Omega \setminus (\Gamma_{\mathrm{1}} \cup \Gamma_{\mathrm{IN}} \cup \Gamma_{\mathrm{OUT}}) $.
The snapshot in figure~\ref{fig:lidriven_snapshot} represents the magnitude of the velocity at time $t=5$s in for the parameter $\lambda=0.5$ as defined in (\ref{eq:lambda_lid_driven}).
}
\end{figure}

Here, we consider the planar flow of an incompressible newtonian fluid in a lid-driven cavity with an inflow and an outflow (see description in figure~\ref{fig:liddriven_configuration}).
We define two Reynolds numbers associated with this configuration.
The first is associated with the lid velocity $\mathrm R _1= {\frac{\rho \, v_{\mathrm{lid}} \, \ell}{\eta}}$ where $\rho=1$kg.m$^{-3}$ is the fluid density, $v_{\mathrm{lid}}=1$m.s$^{-1}$ is the magnitude of the lid velocity $\mathbf{v}\vert_{\Gamma_1}=(1, 0)$, $\ell=1$m is the length of the cavity's sides and $\eta=10^{-3}$kg.m$^{-1}$.s$^{-1}$ is the dynamic viscosity so that $\mathrm{R}_1 = 10^3$.
The second is the Reynolds number associated with the velocity at inflow: $\mathrm R _2= {\frac{\rho \, v_{\mathrm{inflow}} \, \ell}{\eta}}$ with $v_{\mathrm{inflow}}$ the magnitude of inflow velocity $\mathbf{v}\vert_{\Gamma_{\mathrm{IN}}} = (v_x, 0)$.
The dimensionless parameter is the ratio of these two Reynolds numbers 
\begin{equation}
\label{eq:lambda_lid_driven}
\lambda = \frac{\mathrm{R}_2}{\mathrm{R}_1},
\end{equation}
and controls directly the inflow velocity $v_x$.

The HDM is derived exactly as in the previous application (subsection~\ref{sec:canal}).
The mesh includes here 13656 nodes.
Remark that for $\lambda=0$ ($\mathbf{v}\vert_{\Gamma_{\mathrm{IN}}}=(0,0)$), this application corresponds to the standard lid-driven cavity and that the main vortex develops at another location for $0<\lambda\leq 1$.
The test cases 
are given in table~\ref{tab:cases_liddriven}.
The test procedure is identical to that in previous section~\ref{sec:canal}, with the difference that the snapshots are taken from the initial condition over a \emph{transient period} of 10s.
Then, we compute sets of POD bases for the fluctuating velocity and pressure, and we derive interpolated bases by the proposed method and that proposed in~\cite{Amsallem2008}.
Every bases are truncated to $M=20$ modes.

\begin{table}[h!]
	\centering
	\begin{tabular}{ccc}
		\hline
		Case & Sampling $\Lambda$ & Target $\lambda^\star$ \\
		\hline
		1 & (0.0, 0.2, 0.4, 0.6, 0.8, 1.0) & 0.5 \\ 
		2 & (0.0, 0.2, 0.3, 0.4, 0.9, 1.0) & 0.5 \\ 
		3 & (0.0, 0.2, 0.3, 0.4, 0.9, 1.0) & 0.7 \\ 
		4 & (0.0, 0.3, 0.4, 0.9, 1.0) & 0.1 \\ 
		\hline	
	\end{tabular}
	\caption{\label{tab:cases_liddriven} 
	Description of the cases for the application of section~\ref{sec:liddriven} (lid-driven cavity with inflow). The parameter is the ratio of the Reynolds number associated with the lid velocity and that associated with the inflow velocity as defined in (\ref{eq:lambda_lid_driven}). The set of sampling parameters for which the simulation of the full order model is performed is $\Lambda$. The interpolation is performed for the new parameter~$\lambda^\star$.}
\end{table}

We firstly examine the projection error for the fluctuating velocity defined in (\ref{eq:projection_error_definition}). 
The results are given in table~\ref{tab:reconstruction_error_liddriven}, where we see that the error due to the projection of the velocity on the basis interpolated by the proposed method is globally lower than the error due to the projection on the basis interpolated by the method proposed in \cite{Amsallem2008}.
\begin{table}[h!]
\centering
\begin{tabular}{lllll}
\hline
Method & Case 1 & Case 2 & Case 3 & Case 4 \\
\hline
\textsc{Reference} &  2.98e-06  &  2.98e-06  &  4.71e-05  &  1.90e-06  \\
\textsc{Neville} & { 2.7e-02 } & { 2.2e-01 } & { 2.0e-01 } & { 1.6e-01 } \\
\textsc{Amsallem} & { 2.8e-02 } & { 2.1e-01 } & { 4.6e-01 } & { 4.8e-01} \\
\hline
\end{tabular}
\caption{
\label{tab:reconstruction_error_liddriven}
Results for the projection error as defined in (\ref{eq:projection_error_definition}) for the application of section~\ref{sec:liddriven} (lid-driven cavity with inflow).
	\textsc{Reference} refers to the basis built directly from the snapshots obtained by the simulation of the HDM at the target parameter $\lambda ^\star$ (no interpolation). \textsc{Neville} and \textsc{Amsallem} refer to the bases obtained from the interpolation of the sampling sets defined in table~\ref{tab:cases_liddriven} (no simulation of the HDM).
}
\end{table}
Then, we build POD-ROMs from the reference basis \textsc{Reference} and the interpolated bases \textsc{Neville} and \textsc{Amsallem} as described in previous section~\ref{sec:canal}.
The simulation of these reduce-order dynamical systems yields the temporal coefficients $\big(a_i(t)\big)_{i=1}^M$ and $\big(b_i(t)\big)_{i=1}^M$ in (\ref{eq:ersatz_velocity}--\ref{eq:ersatz_pressure}).
%
%
%
The dynamical systems errors defined in (\ref{eq:dynamic_error_definition}) is shown in table~\ref{tab:dynamic_error_liddriven} for the four cases described in table~\ref{tab:cases_liddriven}. We see that the proposed method yields the lowest errors.
\\

\begin{table}[h!]
\centering
\begin{tabular}{lllll}
\hline
Method & Case 1 & Case 2 & Case 3 & Case 4 \\
\hline
\textsc{Reference} &  1.56e-03  &  1.56e-03  &  4.54e-03  &  5.13e-04  \\
\textsc{Neville} & { 4.6e-02 } & { 3.5e-01 } & { 5.1e-01 } & { 2.7e-01 } \\
\textsc{Amsallem} & { 5.1e-02 } & { 3.5e-01 } & { 9.7e-01 } & { 6.2e-01 } \\
\hline
\end{tabular}
\caption{\label{tab:dynamic_error_liddriven}
	Results for the error associated with the dynamical systems as defined in (\ref{eq:dynamic_error_definition}) for the application of section~\ref{sec:liddriven} (lid-driven cavity with inflow).}
\end{table}

\begin{figure}
\centering
\includegraphics[width=0.8\linewidth]{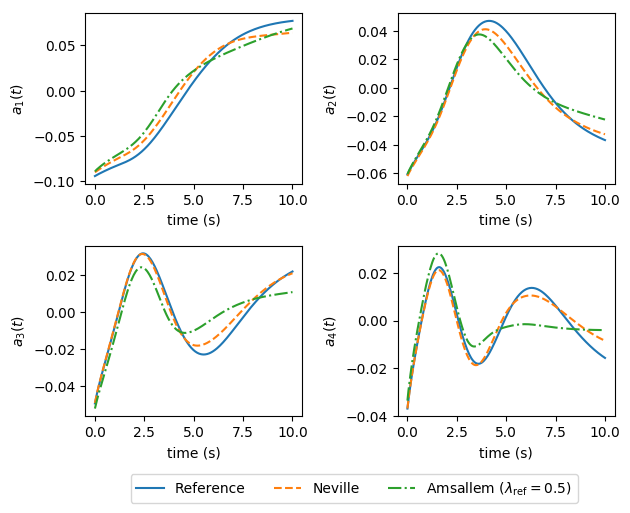} 
\caption{\label{fig:influence_point_de_ref}
Temporal coefficients $\big(a_i(t)\big)_{i=1}^4$ in (\ref{eq:ersatz_velocity}) returned by the dynamical system build from the POD bases obtained directly from the snapshots (\textsc{Reference}) and from the two interpolation methods (\textsc{Neville} and \textsc{Amsallem}) for the application of section \ref{sec:liddriven} (lid-driven cavity with inflow). The sampling parameters are $\Lambda=\{0, 0.2, 0.5, 0.7, 1\}$ and the target parameter is $\lambda^\star = 0.1$. For the method \textsc{Amsallem}, the reference point is chosen as the center of the sample set $\lambda_r = 0.5$.
}
\end{figure}

Finally, we show the result of a bad choice for the reference point $\lambda_r$ in the method from \cite{Amsallem2008} in figure~\ref{fig:influence_point_de_ref}. This also illustrates the benefit of the independence of the proposed method from any reference point.

\subsection{Flows induced by a rotating body\label{sec:rotROM}}
\begin{figure}[h!]
    \centering
    {\hfill
    \begin{subfigure}{.3\linewidth}
        \includegraphics[width=\textwidth]{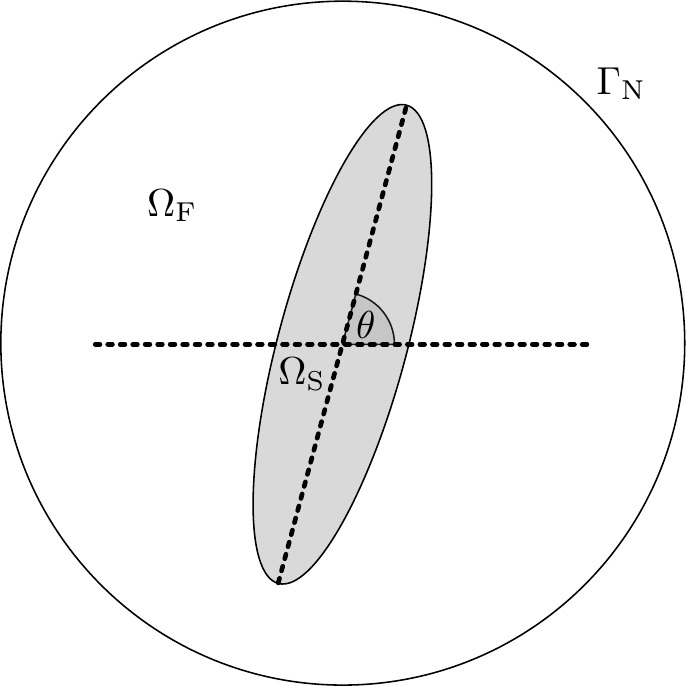}
      
        \caption{\label{fig:rotROM_schematic} Schematic.}
    \end{subfigure}
    \hfill
    \begin{subfigure}{.3\linewidth}
        \includegraphics[width=\textwidth]{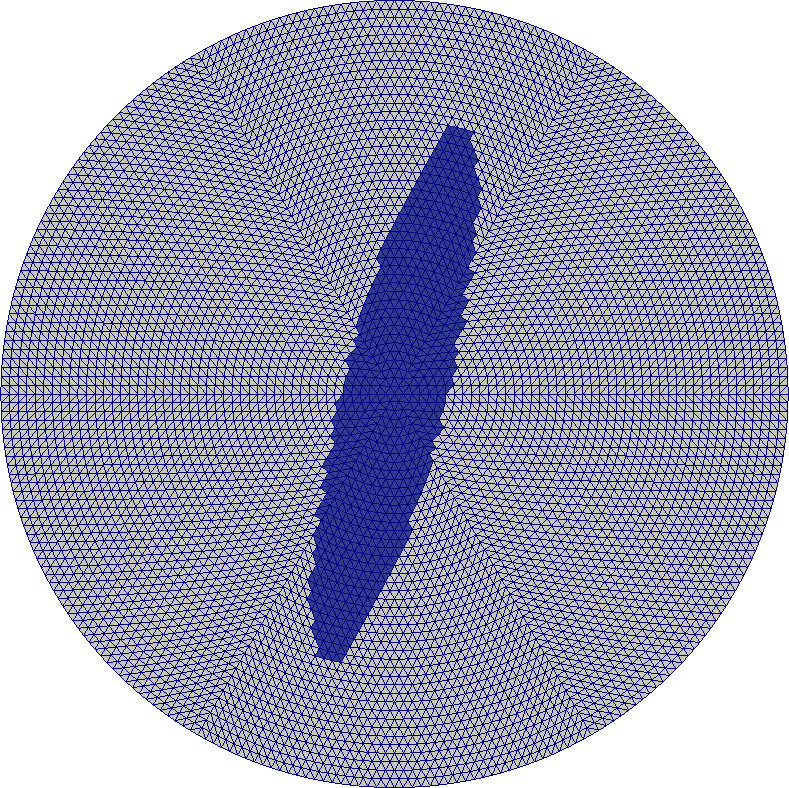}
        \caption{\label{fig:rotROM_mesh} Mesh.}
    \end{subfigure}
    \hfill}    
    \caption{\label{fig:rotROM_configuration} 
    Description of the configuration for the application of section~\ref{sec:rotROM}
     (flow induced by a rotating body). 
The domain radius is $1$m and the ellipse principal radius is $R=0.8$m with an aspect ratio of $0.2$.
The boundary condition is outflow $\boldsymbol{\sigma} \cdot \mathbf{n}=0$ on $\Gamma_{\mathrm{N}}$ (with $\boldsymbol{\sigma}$ the Cauchy stress tensor and $\mathbf{n}$ the outward normal unit vector on the boundary).
%
}
\end{figure}

Here, we consider a two dimensional circular spatial domain $\Omega=\Omega_{\mathrm S} \cup \Omega_{\mathrm F}$ filled with a rotating ellipsoidal body $\Omega_{\mathrm S}$ immersed in an incompressible newtonian fluid $\Omega_{\mathrm F}$ (see description in figure~\ref{fig:rotROM_configuration}).
The parameter is the Reynolds number defined as 
\begin{equation}
\label{eq:reynolds_rotROM}
\lambda = \frac{\rho\,v_\infty\, \ell }{\eta},
\end{equation}
with the density $\rho=1$ (kg.m$^{-3}$), the dynamic viscosity $ \eta=0.01$ (kg.m$^{-1}$.s$^{-1}$), the velocity at the ellipse tips ${v_\infty=R\,\dot\theta}$ (m.s$^{-1}$) and $\ell=2R$ (m) the ellipse principal diameter.
The governing equations are derived by extending the incompressible Navier-Stokes equations to the solid domain by the \emph{fictitious domains} method.
Then, the momentum equation and the continuity equation are solved together by a monolithic formulation on the mixed finite element space known as the \emph{mini space} (\textit{i.e.} linear vector Lagrange element enriched with cubic vector bubble element for velocity and piecewise linear element for pressure, see \cite{arnold1984stable} for details).
The finite-elements mesh includes $52669$ nodes and is not conforming with the body's boundaries (see figure~\ref{fig:rotROM_mesh}).
The test cases (\textit{i.e.} the set of sampling  parameters and the target parameter) are given in table \ref{tab:cases_rotROM}.

\begin{table}[H]
\centering
\begin{tabular}{ccc}
\hline
Case & Sampling $\Lambda$ & Target $\lambda^\star$ \\
\hline
$1$ & $(500, 1000, 1500, 2500, 3000)$ & $2800$ \\ 
$2$ & $(500, 1000, 1500, 1800, 2500, 2800, 3000)$ & $800$ \\ 
$3$ & $(500, 1500, 3000)$ & $2800$ \\
\hline
\end{tabular}
\caption{
\label{tab:cases_rotROM} 
Description of the cases for the application of section~\ref{sec:rotROM} (flow induced by a rotating body). 
The parameter is the Reynolds number as defined in (\ref{eq:reynolds_rotROM}).
The set of sampling parameters for which the simulation of the full order model is performed is $\Lambda$. 
The interpolation is performed for the new parameter~$\lambda^\star$.
}
\end{table}

The test procedure is identical to that described in section~\ref{sec:canal}, except for the construction of the POD-ROM which has to comply with the multiphase description of the fluid--structure interaction (the interested reader is referred to~\cite{falaize2018pod}).
As in previous section~\ref{sec:liddriven}, we are concerned here with the \emph{transient period} (half rotation of the ellipse).
We compute sets of POD bases for the fluctuating velocity, and we derive interpolated bases by the proposed method and that proposed in~\cite{Amsallem2008}.
Every bases are truncated here to $M=40$ modes.

\begin{table}[h!]
\centering
\begin{tabular}{llll}
\hline
Method & Case 1 & Case 2 & Case 3\\
\hline
\textsc{Reference} &  6.493e-04  &  2.747e-04  &  6.493e-04   \\
\textsc{Neville} & { 1.874e-03 } & { 2.457e-03 } & { 2.141e-03 } \\
\textsc{Amsallem} &  { 1.880e-03 } & { 2.457e-03 } & { 2.143e-03 } \\
\hline
\end{tabular}
\caption{
\label{tab:reconstruction_error_rotROM}
Results for the projection error as defined in (\ref{eq:projection_error_definition}) for the application of section~\ref{sec:rotROM} (flow induced by a rotating body).
}
\end{table}

\begin{figure}[H]
\centering
\includegraphics[width=0.7\linewidth]{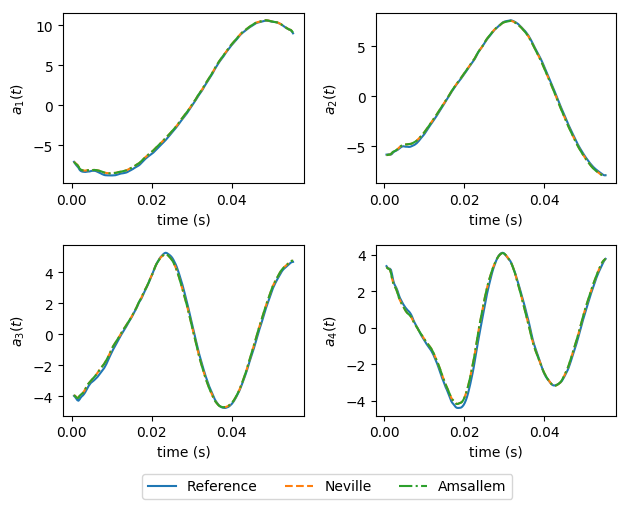}
\caption{
		\label{fig:rotROM_coeffs} 
		Temporal coefficients $\big(a_i(t)\big)_{i=1}^4$ in (\ref{eq:ersatz_velocity}) returned by the dynamical system build from the reference POD basis obtained directly from the snapshots (\textsc{Reference}) and from the two interpolation methods (\textsc{Neville} and \textsc{Amsallem}) for the application of section \ref{sec:rotROM} (flow induced by a rotating body) in the case~2 (see table~\ref{tab:cases_rotROM}).
}
\end{figure}

The projection error for the fluctuating velocity defined in (\ref{eq:projection_error_definition}) are given in table~\ref{tab:reconstruction_error_rotROM}.
The error associated with the reduced-order dynamical systems (in table~\ref{tab:dynamic_error_rotROM}).
Globally, the proposed method yields the lowest errors.
Additionally, the temporal coefficients $\big(a_i(t)\big)_{i=1}^M$ for the fluctuating velocity are shown in figure~\ref{fig:rotROM_coeffs}.

\begin{table}[h!]
\centering
\begin{tabular}{llllll}
\hline
Method & Case 1 & Case 2 & Case 3 \\
\hline
\textsc{Reference} &  2.964e-02  &  3.029e-03  &  2.964e-02   \\
\textsc{Neville} & { 2.108e-02 } & { 1.817e-02 } & { 1.474e-02 } \\
\textsc{Amsallem} &  {2.189e-02 } & { 1.818e-02 } & { 1.475e-02 } \\
\hline
\end{tabular}
\caption{
\label{tab:dynamic_error_rotROM}
Results for the error associated with the dynamical systems as defined in (\ref{eq:dynamic_error_definition}) for the application of section~\ref{sec:rotROM} (flow induced by a rotating body).
	}
\end{table}

%

\section{Conclusion} In this work, we extend the well-known Neville-Aitken's interpolation algorithm to 
the Grassmann manifold for the construction of reduced order methods. The recursive character of the algorithm allows us to perform interpolation on the Grassmann manifold using the geodesic barycenter of two points.

Our method presents several advantages with respect to the algorithms developed in \cite{Amsallem2008,Mosquera_DCDS18}. Indeed, it does not require a reference point as in \cite{Amsallem2008}. Moreover, the proposed method is direct and does not require the resolution of a fixed point problem for the minimization process as in \cite{Mosquera_DCDS18} in the framework of the IDW method on Grassmann manifolds. Furthermore, our algorithm is more 
pertinent both in terms of the accuracy of results and the time computation. The performances of our proposed method are illustrated through three independent CFD applications, namely: the Von Karman vortex shedding street, the lid-driven cavity with inflow and the flow induced by a rotating solid.


\end{document}